\documentclass[12pt]{amsart}
\usepackage{amsmath, amsthm, amsfonts, amssymb}

\newtheorem{theorem}{Theorem}[section]
\newtheorem{lemma}[theorem]{Lemma}

\newtheorem{corollary}[theorem]{Corollary}
\theoremstyle{definition}

\newtheorem{example}[theorem]{Example}

\theoremstyle{remark}
\newtheorem{remark}[theorem]{Remark}
\numberwithin{equation}{section}

\newcommand{\mC}{\ensuremath{\mathbb{C}}}
\newcommand{\mD}{\ensuremath{\mathbb{D}}}

\newcommand{\mN}{\ensuremath{\mathbb{N}}}

 \allowdisplaybreaks
 
\begin{document}

\title{A normality criterion generalizing Gu's result}

\author[K. S. Charak]{ Kuldeep Singh Charak}
\address{Department of Mathematics, University of Jammu,
Jammu-180 006, India}
\email{kscharak7@rediffmail.com}

\author[V. Singh]{Virender Singh}
\address{Department of Mathematics, University of Jammu,
Jammu-180 006, India}
\email{virendersingh2323@gmail.com }

\begin{abstract} In this paper we prove a normality criterion for the families of meromorphic functions involving sharing of functions. Our result generalizes some of the earlier results on Gu's normality criterion.
\end{abstract}

\renewcommand{\thefootnote}{\fnsymbol{footnote}}
\footnotetext{The research work of the second author is supported by the CSIR India.}

\maketitle

\section{\textbf{Introduction and Main Results}}
It is assumed that the reader is familiar with the standard
notions used in the Nevanlinna value distribution theory such as
$T(r,f),m(r,f),N(r,f),S(r,f)$ etc., one may refer
to \cite{hayman-1}.

\medskip

A family $\mathcal F$ of meromorphic functions defined on a domain $D \subseteq \overline {\mC}$ is said to be normal in $D$ if every sequence of elements of $\mathcal F$  contains a subsequence which converges locally uniformly in $D$ with respect to the spherical metric, to a meromorphic function or $\infty$ (see \cite{schiff-1}). 

\medskip

 Two nonconstant meromorphic functions $f$ and $g$ defined on $D$ are said to share a meromorphic functions $\psi$ in $D$ if $\overline{E}_f(\psi)=\overline{E}_g(\psi)$, where
$$\overline{E}_f(\psi)=\{z\in D: f(z)=\psi(z)\}.$$

\medskip

The following Picard type theorem is one of the main result from Hayman's seminal paper \cite{hayman-2}:
 \begin{theorem}(Hayman's alternative) Let $f$ be a nonconstant meromorphic function in $\mC$, $k$ a natural number and $c$ a nonzero complex number. Then $f$ or $f^{(k)}-c$ has a zero in $\mC$. If $f$ is transcendental, $f$ and $f^{(k)}-c$ has infinitely many zeros in $\mC$.
\end{theorem}
In 1979, Y.X. Gu \cite{gu-1} proved the following normality criterion corresponding to Hayman's alternative:
\begin{theorem}\label{theorem1}(Gu's normality criterion) Let $\mathcal{F}$ be a family of meromorphic functions defined in a domain $D$, and let $k$ be a positive integer. If, for every function $f\in\mathcal{F}$, $f\neq 0$ and $f^{(k)}\neq 1$ in $D$, then $F$ is normal in $D$.
\end{theorem}
Since then many variations of Theorem \ref{theorem1} have been obtained, for instance one can see \cite{fang-2, liu-1, nevo-1, schwick-2, yang-1,yang-2}.
In fact Schwick \cite{schwick-1} proved a more general version of Gu's result:

\begin{theorem}\label{theorem2} Let $\psi\not\equiv 0$ be a meromorphic function in a domain $D$ and $k \in \mN $. Let $\mathcal{F}$ be a family of meromorphic functions in $D$, such that $f\neq 0$ and $f^{(k)}\neq \psi $, and $f$ and $\psi$ have no common poles for each $f \in \mathcal{F}$. Then $\mathcal{F}$ is normal in $D$.
\end{theorem}

In 2004, Fang and Zalcman \cite{fang-1} proved the following generalization of Theorem \ref{theorem1} by considering the sharing of values:
\begin{theorem} \label{theorem3}
Let $k$ be a positive integer and $b$ be a nonzero complex constant. Let $\mathcal{F}$ be a family of meromorphic functions on $D$, all of whose zeros have multiplicity at least $k+2$, such that for each pair of functions $f$ and $g$ in $\mathcal{F}$, $f$ and $g$ share the value $0$, and $f^{(k)}$ and $g^{(k)}$ share the value $b$ in $D$, then $\mathcal{F}$ is normal in $D$.
\end{theorem}

Recently, J. Chang \cite{chang-1} proved the following result by replacing the constant $b$ by a holomorphic function:

\begin{theorem} \label{theorem4} Let $k \in \mN$ and $h(\not\equiv 0)$ be a function holomorphic on $D$. Let $\mathcal{F}$ be a family of meromorphic functions in $D$, all of whose zeros have multiplicity at least $k+2$, such that for each pair of functions $f$ and $g$ in $\mathcal{F}$, $f$ and $g$ share the value $0$, and $f^{(k)}$ and $g^{(k)}$ share the function $h$. Suppose additionally that at each common zero of $f$ and $h$ for every $f \in \mathcal{F}$, the multiplicities $m_f$ for $f$ and $m_h$ for $h$ satisfy $m_f\geq m_h +k+1$ for $k>1$ and $m_f \geq 2m_h +3$ for $k=1$. Then, $F$ is normal in $D$.
\end{theorem}

Examples are also given in \cite{chang-1} for the sharpness of conditions in Theorem \ref{theorem4}. Working in this direction, we prove the following generalization of Theorem \ref {theorem4} :
\begin{theorem}\label{theorem5} Let $\mathcal{F}$ be a family of meromorphic functions in a domain $D$, and let $k$ be a positive integer. Suppose that $\phi$ is a holomorphic function on $D$ and $\psi$ is a meromorphic function on $D$ such that $\phi^{(k)} (z)\not\equiv \psi (z)$. Suppose that for each pair of functions $f$ and $g$ in $\mathcal{F}$, $f$ and $g$ share $\phi$, and $f^{(k)}$ and $g^{(k)}$ share the function $\psi$.
Suppose further that
\begin{enumerate}
\item every $f \in\mathcal{F}$, $f-\phi$ has zeros of multiplicity at least $k+2$,
\item for every common zero of $f-\phi$ and $\psi -\phi^{(k)}$, the multiplicities $m_{f-\phi}$ for $f-\phi$ and $m_{\psi-\phi^{(k)}}$ for $\psi -\phi^{(k)}$ satisfy $$m_{f-\phi} \geq m_{\psi-\phi^{(k)}}+k+1 \text{~for~} k>1,\text{~and}$$ $$m_{f-\phi} \geq 2 m_{\psi-\phi^{(k)}}+3 \text{~for~} k=1, $$
\item for every $f \in \mathcal{F}$, $f$ and $\psi$ have no common poles in $D$.
\end{enumerate}
Then $\mathcal{F}$ is normal on $D$.
\end{theorem}

\begin{remark} If $\phi\equiv 0 $ and $\psi$ is a holomorphic function, then Theorem \ref{theorem5} reduces to Theorem \ref{theorem4}. Thus, the conditions (1) and (2) in Theorem \ref{theorem5} can easily be seen to be essential.
\end{remark}

\begin{example}Consider the family $$\mathcal F =\left\{\frac{1}{2mz}:m\in \mN\right\}$$ on the open unit disk $\mD$, and let $\phi(z)=1/z$ and $\psi (z)\equiv 0$. Then clearly, for every $f,g \in \mathcal F$, $f$ and $g$ share $\phi(z)$, and $f^{(k)}$ and $g^{(k)}$ share $\psi(z)$ in $\mD$. However, the family $F$ is not normal in $\mD$. This shows that $\phi$ cannot be taken meromorphic in Theorem \ref{theorem5}. 

\medskip

Further, for the same family $\mathcal F$, if we take $\phi(z)\equiv 0$ and $\psi(z)=1/z^{k+1}$, then for every $f,g\in \mathcal F$, $f$ and $g$ share $\phi(z)$, and $f^{(k)}$ and $g^{(k)}$ share $\psi(z)$ in $\mD$. But, the family $\mathcal F$ is not normal in $\mD$. This shows that the condition (3) in Theorem \ref{theorem5} is essential.
\end{example}

\begin{example}Consider the family $$\mathcal F =\left\{z^{k+1}+\frac{1}{mz}:m\in \mN\right\}$$ on the open unit disk $\mD$, and let $\phi(z)=z^{k+1}$ and $\psi (z)=(k+1)!z$. Then clearly, for every $f,g \in \mathcal F$, $f$ and $g$ share $\phi(z)$, and $f^{(k)}$ and $g^{(k)}$ share $\psi(z)$ in $\mD$. However, the family $F$ is not normal in $\mD$. This shows that the condition $\phi^{(k)}(z)\not\equiv \psi(z)$ in Theorem \ref{theorem5} cannot be dropped. 
\end{example}

\begin{corollary}\label{corollary1} Let $n \geq 2$ be a positive integer and $\psi(\not\equiv 0)$ be a function meromorphic on $D$. Let $\mathcal{F}$ be a family of meromorphic functions in $D$ such that for each pair of functions $f$ and $g$ in $\mathcal{F}$, $f$ and $g$ share the value $0$, and $f^n f'$ and $g^n g'$ share the function $\psi$. Suppose further
\begin{enumerate}
\item for every common zero of $f$ and $\psi$, the multiplicities $m_f$ for $f$ and $m_{\psi}$ for $\psi$ satisfy $$m_{f} \geq m_{\psi}+k+1 \text{~for~} k>1,\text{~and}$$ $$m_{f} \geq 2 m_{\psi}+3 \text{~for~} k=1, $$
\item for every $f \in \mathcal{F}$, $f$ and $\psi$ have no common poles in $D$.
\end{enumerate}
Then $\mathcal{F}$ is normal on $D$.
\end{corollary}

\medskip

Corollary \ref{corollary1} follows by setting $\mathcal {G}=\left\{f^{n+1}/(n+1): f\in \mathcal {F}\right\}$ and applying Theorem \ref{theorem5} to this family with $\phi(z)\equiv 0$ and $k=1$.

\section{Proof of the Main result}
For $z_0\in\mathbb{C}$ and $r>0$, we denote by $D_r(z_0)$ the open unit disk with centre $z_0$ and radius $r$, and $D' _r(z_0)$ the corresponding punctured disk. To prove our main result-Theorem \ref{theorem5}, we require the following lemma:

\begin{lemma}\label{lemma1}\cite{chen-1} Let $\mathcal{F}$ be a family of meromorphic functions in a domain $D$, all of whose zeros have multiplicity at least $k$. Then, if $\mathcal{F}$ is not normal at $z_0$, then for each $\beta:-1 < \beta <k$, there exist points $z_n \in D$ with $z_n\rightarrow z_0$, functions $f_n \in \mathcal{F}$ and positive numbers $\rho_n \rightarrow 0$ such that $$g_n(\zeta):=\rho^{-\beta}_n f_n (z_n+ \rho_n \zeta)$$ converges locally uniformly with respect to the spherical metric in $\mC$ to a nonconstant meromorphic function $g$ of finite order, all of whose zeros have multiplicity at least $k$.
\end{lemma}

Further, recall that the sets $\left\{E_{\lambda}\right\}_{\lambda \in \Lambda}$ are said to be {\it locally uniformly discrete} in $D$, if for each point $z_0 \in D$, there exists $\delta>0$ such that either $E_{\lambda}\cap D_\delta(z_0)$ is an empty set or a singleton, one may refer to \cite{chang-1}.

\begin{proof} [\textbf{Proof of Theorem \ref{theorem5}.}]  Since normality is a local property, it is enough to show that $\mathcal{F}$ is normal at each $z_0\in D$. We distinguish the following cases.\\\\
\textbf{Case I.} Suppose that there exist $f\in \mathcal{F}$ such that $f(z_0)\neq \phi(z_0)$ and $f^{(k)}(z_0)\neq \psi(z_0)$. Then we can find $r>0$ such that $D_r(z_0)\subset D$, and $f(z)\neq \phi(z) $ and $f^{(k)}(z)\neq \psi(z)$ in $D_r(z_0)$ and so by the given sharing condition, $f\neq \phi$ and $f^{(k)} \neq \psi$ for every $f\in \mathcal F$ in $D_r(z_0)$. Now set $\alpha:=\psi - \phi^{(k)}$ and consider the family $$\mathcal G=\left\{g=f-\phi:f\in \mathcal{F}\right\}.$$ Then clearly $\alpha(z)\not\equiv 0$ and for every $g\in \mathcal G$, $g(z)\neq 0$ and $g^{(k)}(z)\neq \alpha (z)$. Thus by Theorem \ref{theorem2} $\mathcal G$ is normal in $D_r(z_0)$. Since $\mathcal G$ is normal if and only if $\mathcal F$ is normal, $\mathcal F$ is normal at $z_0$.\\\\
\textbf{Case II.} Suppose that there exist $f\in \mathcal F$ such that $f(z_0)=\phi(z_0)$ or $f^{(k)}(z_0)=\psi(z_0)$. Then we can find $r>0$ such that $D_r(z_0)\subset D$, and $f(z)\neq \phi(z)$ and $f^{(k)}(z)\neq \psi(z)$ in $D'_r(z_0)$ and so by the given sharing condition, $f(z)\neq \phi(z)$ and $f^{(k)}(z)\neq \psi(z)$ for every $f\in \mathcal F$ in $D'_r(z_0)$. Thus for any $z_1 \in D_r(z_0)$, there exists $\delta>0$ such that every $E_f$ has at most one point lying in $D_\delta (z_1)$, where $$E_f:=\left\{z\in D_r(z_0):f(z)=\phi(z)\right\}\cup \left\{z\in D_r(z_0):f^{(k)}(z)=\psi(z)\right\}.$$ Therefore, the sets $\left\{E_f\right\}_{f\in \mathcal F}$ are locally uniformly discrete in $D_r(z_0)$.

\medskip

 As in Case I, consider the family $\mathcal G=\left\{g=f-\phi:f\in \mathcal{F}\right\}$ and set $\alpha:=\psi - \phi^{(k)}$. Then clearly the sets $\left\{E_g\right\}_{g\in \mathcal G}$ are locally uniformly discrete in $D_r(z_0)$, where $$E_g=\left\{z\in D_r(z_0):g(z)=0\right\}\cup \left\{z\in D_r(z_0):g^{(k)}(z)=\alpha(z)\right\}.$$ If $\alpha(z)$ is holomorphic in $D_r(z_0)$, then by \cite[Theorem 4, p.49]{chang-1}, $\mathcal G$ is normal in $D_r(z_0)$ and hence $\mathcal F$ is normal at $z_0$. Suppose that $\alpha(z)$ is not holomorphic in $D_r(z_0)$. Assume that $z_0$ is a pole of $\alpha(z)$. Then we can find $\delta>0$ such that $D_\delta (z_0) \subset D_r(z_0)$ and $\alpha(z)$ is holomorphic in $D'_\delta(z_0)$, and thus $\mathcal G$ is normal in $D'_\delta(z_0)$. Next, consider the family $$\mathcal H:=\left\{h(z)=\frac{g(z)}{\alpha(z)}:g\in \mathcal G\right\}.$$ Noting that $z_0$ is a pole of $\alpha(z)$ and for every $f\in \mathcal F$, $f$ and $\psi$ have no common poles implies that\\
$(a)$ for every $g \in \mathcal G$, $g$ and $\alpha$ have no common poles and hence for every $h\in \mathcal H$, $z_0$ is a zero of $h$ of multiplicity at least $k+3$,\\
$(b)$ there exists $\eta >0$ such that $D_{\eta}(z_0)\subset D_\delta(z_0)$ and for every $h \in \mathcal H$, $h\neq 1,\infty$ in $D_{\eta}(z_0)$.

\medskip

We first prove that $\mathcal H$ is normal at $z_0$. Suppose on contrary that $\mathcal H$ is not normal at $z_0$. Then by Lemma \ref{lemma1}, we can find a sequence $\left\{ h_j \right\}$ in $\mathcal{H}$, a sequence $\left\{ z_j\right\}$ of complex numbers with $z_j\rightarrow 0$ and a sequence $\left\{\rho_j\right\}$ of positive real numbers with $\rho_j \rightarrow 0 $ such that $$H_j (\zeta) = h_j (z_j +\rho_j \zeta)$$ converges locally uniformly with respect to the spherical metric to a nonconstant meromorphic function $H(\zeta)$ on $\mC$, all of whose zeros have multiplicity at least $k+3$. Also by Hurwitz theorem, we have $H\neq 1,\infty$ on $\mC$. Thus by second fundamental theorem of Nevanlinna, we have
\begin{align*} T(r,H)&\leq \overline N(r,H) + \overline N\left(r,\frac{1}{H}\right)+ \overline{N}\left(r,\frac{1}{H-1}\right)+S(r,H)\\
                     &\leq \frac{1}{k+3}N\left(r,\frac{1}{H}\right) + S(r,H)\\
										 & \leq \frac{1}{k+3}T(r,H)+ S(r,H),
\end{align*}
which is a contradiction. Therefore $\mathcal H$ is normal at $z_0$. Now we turn to prove the normality of $\mathcal G$ at $z_0$.

\medskip

 Suppose that $\mathcal{G}$ is not normal at $z_0$. Since $\mathcal H$ is normal at $z_0$, it is equicontinuous at $z_0$ with respect to the spherical metric. Also $h(z_0)=0$ for every $h\in \mathcal H$. Thus there exists $\delta_1>0$ such that $D_{\delta_1}(z_0)\subset D_{\delta}(z_0)$ and $|h(z)|\leq 1$ for every $h\in \mathcal H$ in $D_{\delta_1}(z_0)$. It follows that $\mathcal G$ is a family of holomorphic functions in $D_{\delta_1}(z_0)$. 

\medskip

Let $\left\{g_n\right\}$ be a sequence in $\mathcal G$. Since $\mathcal{G}$ is normal in $D'_{\delta_1} (z_0)$ but not at $z_0$, there exists a subsequence of $\left\{g_n\right\}$, which we may take as $\left\{g_n\right\}$ itself, which converges locally uniformly on $D'_{\delta_1} (z_0)$ but not on $D_{\delta_1} (z_0)$. By the maximum modulus principle, we have $\left\{g_n\right\}$ converges locally uniformly to $\infty$ in $D'_{\delta_1}(z_0)$ and hence $\left\{h_n\right\}$ converges locally uniformly to $\infty$ in $D'_{\delta_1}(z_0)$, which is a contradiction to the fact that $|h(z)|\leq 1$ for every $h\in \mathcal H$ in $D_{\delta_1}(z_0)$.
Thus $\mathcal G$ is normal at $z_0$ and hence $\mathcal F$ is normal at $z_0$. 

\end{proof}

\bibliographystyle{amsplain}

\end {document}